\documentclass[anon,12pt]{colt}

\usepackage{graphicx} 
\usepackage{array}
\usepackage{paralist}

\usepackage{amsmath,amssymb,latexsym}
\usepackage{algorithm,algorithmic}

\newcommand{\BEAS}{\begin{eqnarray*}}
\newcommand{\EEAS}{\end{eqnarray*}}
\newcommand{\BEA}{\begin{eqnarray}}
\newcommand{\EEA}{\end{eqnarray}}
\newcommand{\BEQ}{\begin{equation}}
\newcommand{\EEQ}{\end{equation}}
\newcommand{\BIT}{\begin{itemize}}
\newcommand{\EIT}{\end{itemize}}
\newcommand{\BNUM}{\begin{enumerate}}
\newcommand{\ENUM}{\end{enumerate}}
\newcommand{\QED}{\hspace*{\fill}~$\blacksquare$}

\title[Streaming Matrix completion]{Streaming, Memory Limited Matrix Completion with Noise}
\usepackage{times}

\coltauthor{\Name{Se-Young Yun} \Email{seyoung.yun@inria.fr}\\
\Name{Marc Lelarge} \Email{marc.lelarge@ens.fr}\\
\Name{Alexandre Proutiere} \Email{alepro@kth.se}
}

\begin{document}
\maketitle
\begin{abstract}
In this paper, we consider the streaming memory-limited matrix completion problem
when the observed entries are noisy versions of a small random
fraction of the original entries. We are interested in scenarios where the matrix size
is very large so the matrix is very hard to store and manipulate. Here, columns of the observed matrix
are presented sequentially and the goal is to complete the missing
entries after one pass on the data with limited memory space and
limited computational complexity. We propose a streaming algorithm
which produces an estimate of the original matrix with
a vanishing mean square error,
uses memory space scaling linearly with the ambient dimension of the
matrix, i.e. the memory required to store the output alone, and spends
computations as much as the number of non-zero entries of the input
matrix. 
\end{abstract}

\begin{keywords}
matrix completion,
streaming input,
limited memory,
computational complexity
\end{keywords}

\section{Introduction}

Reconstructing a structured (e.g. low rank) matrix from noisy observations of a subset of its entries constitutes a fundamental problem in collaborative filtering \cite{rennie2005}, and has recently attracted much interest, see e.g. \cite{candes2009exact}, \cite{candes2010power}, \cite{keshavan2010matrix}, \cite{recht2011simpler}. The recent development of matrix completion algorithms has been largely motivated by the design of efficient recommendation systems. These systems (amazon, netflix, google) aim at proposing items or products from large catalogues to targeted users based on the ratings provided by users of a small subset of items. This goal naturally translates to a matrix completion problem where the rows (resp. the columns) of the matrix correspond to items (resp. to users). And often, the (item, user) rating matrix is believed to exhibit a low rank structure due to the inherent similarities among users and among items. 

In this paper, we address the problem of matrix completion in scenarios where the matrix can be extremely large, so that (i) it might become difficult to manipulate or even store, and (ii) the complexity of the proposed algorithms should not rapidly increase with the matrix dimensions. In other words, we aim at designing matrix completion algorithms under memory and computational constraints. Memory-limited algorithms are particularly relevant in the streaming data model, where observations (e.g. ratings in recommendation systems) are collected sequentially. We assume here that the columns of the matrix are revealed one by one to the algorithm. More specifically, a subset of noisy entries of an arriving column is observed, and may be stored, but the algorithm cannot request these entries later if they were not stored. The streaming model seems particularly appropriate to model recommendation systems, where users actually seek for recommendations sequentially. Recently, motivated by the need to understanding high-dimensional data, several machine learning techniques, such as PCA \cite{mitliagkas2013} or low-rank matrix approximation \cite{clarkson2009numerical}, have been revisited considering memory and computational constraints. To our knowledge, this paper provides the first analysis of the matrix completion problem under these constraints (refer to the related work section for a detailed description of the connection of our problem to existing work).

Throughput the paper, we use the following notations. For any $m\times n$ matrix $A$, we denote by $A^\dagger$ its transpose. We also denote by $s_1(A)\geq \dots \geq  s_{n\wedge m}(A)\geq 0$, the singular values of $A$. The SVD of matrix $A$ is $A=U\Sigma V^{\dagger}$ where $U,V$ are unitary matrices and $\Sigma = \mbox{diag}(s_1(A),\dots s_{n\wedge m}(A))$. $A^{-1}$ denotes the Pseudo-inverse matrix of $A$, i.e. $ A^{-1} = V \Sigma^{-1}U^\dagger$. Finally, for any vector $v$, $\| v\|$ denotes its Euclidean norm, whereas for any matrix $A$, $\|A
  \|_F$ denotes its Frobenius norm, $\|A\|_2$ its operator norm, and $\|A\|_\infty$ its $\ell_\infty$-norm, i.e., $\|A\|_\infty=\max_{i,j}|A_{ij}|$.

\medskip
\noindent
{\bf Contributions.} Let $M\in [0,1]^{m\times n}$ denote the $m\times n$ ground-truth matrix we wish to recover from noisy observations of some of its entries. $M$ is assumed to exhibit a {\it sparse structure} (refer to Assumption 1 (i) for a formal definition). $m$ and $n$ are typically very large, and can be thought as tending to $\infty$. The SVD of $M$ is $M=U \Lambda V^{\dagger} $. We assume that each entry of $M$ is observed (but corrupted by noise) with probability $\delta$ (independently over entries). The random set of observed entries is denoted by $\Omega$, and we introduce the following operator from $\mathbb{R}^{m\times n}$ to itself: for all $Y\in\mathbb{R}^{m\times n}$,  
$$
[\mathcal{P}_{\Omega}(Y)]_{ij} = \begin{cases} Y_{ij}, \quad & \mbox{if}\quad
  (i,j)\in \Omega \cr 0,& \mbox{otherwise} .\end{cases} 
  $$
Then, we wish to reconstruct $M$ from the observed matrix $A = \mathcal{P}_{\Omega}(M+X)$, where $X$ is a noise matrix with independent and zero-mean entries, and such that $M_{ij}+X_{ij} \in [0,1]$. Note that $\delta$ typically depends of $n$ and $m$, and tends to zero as $n$ and $m$ tend to infinity. Finally, we analyze the matrix completion problem under the {streaming model:} we assume that in each round, a column of $A$ is observed. This column is uniformly distributed among the set of columns that have not been observed so far. 

We present SMC (Streaming Matrix Completion), a memory-limited and low-complexity algorithm which, based on the observed matrix $A$, constructs an estimator $\hat{M}$ of $M$. We prove, under mild assumptions on $M$ and the proportion $\delta$ of observed entries, that $\hat{M}$ is {\it asymptotically accurate}, in the sense that its average mean-square error converges to 0 as both $n$ and $m$ grows large, i.e., $\frac{\|\hat{M}-M \|_F^2}{mn} = o(1)$. More precisely, we make the following assumption. 

\medskip
\noindent
{\it Assumption 1.} (i) $\|M\|^2_F = \Theta(mn)$. (ii) (Structural sparsity of $M$) there exists $i\leq \min(n,m)$ such that $\frac{s_i(M)}{s_{i+1}(M)}=\omega(1)$ and $\sum_{j=i+1}^{K}s_{j}^2(M) = o(mn)$. We denote by $k$ the smallest $i$ satisfying this condition. (iii) $\delta = \omega (k \max({1\over n},\frac{\log^2 m}{ m}))$, and  $\delta = o(\frac{1}{\log^2 m})$. 

\medskip
\noindent
The main result of this paper is a direct consequence of Theorems~\ref{thm:stmcomp}, \ref{thm:memory}, and \ref{thm:complex}. It states that under {Assumption 1}, with high probability, the SMC algorithm provides an asymptotically accurate estimate $\hat{M}$ of $M$ using one pass on the observed matrix $A$, and requires $O(km+kn)$ memory space and $O(\delta mn k)$ operations.

Note that Assumption 1 (ii) is satisfied as soon as $M$ has low rank. More precisely, when rank$(M)=K$, then (ii) is satisfied when $k=K$. In such a case, there is a non-empty set of sampling rates $\delta$ for which SMC yields an asymptotically accurate estimate of $M$ as soon as $K=o({\sqrt{m}\over \log(m)^2})$ (if for example $m$ and $n$ grows at the same pace to infinity). 



\noindent{\bf \em Algorithm overview:} We fully exploit the sparsity of
the input matrix.
\begin{itemize}
\item At first, we store $\ell$ arriving columns into the
memory space. Because of the sparsity, when $\ell = O(\delta
(m+n)/\log m)$, the required memory space to store $\ell$ columns is
$O(m+n)$. While such $\ell$ is not enough to find principle column
vectors, we can reconstruct $K$ principle row vectors for the $m\times
\ell$ submatrix of $M$ corresponding to the arriving columns
(Theorem~\ref{pro:pca}). 
\item Let $A^{(B)}$ denote the first $\ell$ arriving vectors and
  $Q^\dagger$ be the $K$ principle row vectors reconstructed from
  $A^{(B)}$. Using $Q$ and $A^{(B)}$, we generate $K$ reference column
  vectors $W = A^{(B)} \cdot Q$. Although the columns of $W$ are not
  the principle column vectors of $M$, they are very useful. Indeed,
  $A$ has $K$ strong singular values which are much larger than others
  and stem from $M$. Therefore, $\hat{V} = A^\dagger \cdot W$ amplifies
  the part of $W$ sharing the same column space with $M$ a lot more
  than the remaining part. In the sequel, the rows of
  $\hat{V}^\dagger$ become very close to the principle row vectors of
  $M$ (Theorem~\ref{thm:streamingV}).
\item Once we know $\hat{V}^\dagger$ of which rows are very close to the principle row
  vectors of $M$, it is very easy to find
  column vectors $\hat{U}$ such that
  $\frac{\|\hat{U}\hat{V}-M\|_F^2}{mn}=o(1)$. Using the Gram-Schmidt
  process, we can easily find $\hat{R}$ such that $\hat{V}\hat{R}$ is an orthonormal
  matrix and $\hat{U} = \frac{1}{\delta}A \hat{V} \hat{R} \hat{R}^\dagger$ is the
  solution since $\hat{U}\hat{V}^\dagger = \frac{1}{\delta} A \hat{V}
  \hat{R} (\hat{V} \hat{R})^\dagger$ where $\hat{V}
  \hat{R} (\hat{V} \hat{R})^\dagger$ is the projection matrix onto the
  linear span of the rows of $M$.
\end{itemize}
For this algorithm, 1-pass streaming input is enough. For the
streaming input, we can sequentially compute $\hat{V}^j = A_j^\dagger
W$ and update $\hat{U} =
(\sum_{i=1}^{n} A_i \hat{V}^i )\cdot \hat{R} \hat{R}^\dagger$ where $\hat{R}
\hat{R}^\dagger$ can be computed after we get $\hat{V}$. Therefore, we do
not need to store columns to find $\hat{V}$ and $\hat{U}$ for the
1-pass streaming input except the first $\ell$ arriving columns.

\medskip
\noindent{\bf Additional Notations.} When matrices $A$ and $B$ have the same number of rows, $[A, B]$ to denote the matrix whose first columns are those of $A$ followed by those of $B$. For any matrix $A$, $A_{\bot}$ denotes an orthonormal basis of the subspace perpendicular to the linear span of the columns of $A$. $A_i$, $A^j$, and $A_{ij}$ denote the $i$-th column of $A$, the $j$-th row of $A$, and the $(i,j)$ entry of $A$, respectively. For $b\ge a$, $A^{a:b}$ and $A_{a:b}$ are submatrices of $A$ respectively defined as $A^{a:b}=(A^j)_{j=a,\ldots,b}$ and $A_{a:b}=(A_i)_{i=a,\ldots,b}$. Finally, we define the following thresholding operator for matrices. The operator is defined by two real positive numbers $a$ and $b$, with $b\ge a$, and if applied to $A$, it returns the matrix $|A|_a^b$ such that 
$$[|A|^b_a]_{ij} = \begin{cases} b \quad & \mbox{if}\quad
  A_{ij} \ge b ,\cr A_{ij}\quad & \mbox{if}\quad a<A_{ij}<b,\cr
a \quad & \mbox{if}\quad A_{ij}\le a.\end{cases} $$

\section{Related Work}

This section surveys existing work on the design of matrix completion algorithms. We also provide a description of recent work on rank-$k$ approximation and PCA algorithms, as these algorithms could be seen as building blocks of matrix completion methods. The section is organised as follows. We first review algorithms for matrix completion. We then focus on streaming algorithms for rank-$k$ approximation, and PCA. Finally we discuss algorithms designed to be computationally efficient. 

\medskip
\noindent
{\bf Matrix completion algorithms.} \cite{candes2009exact} first showed that in absence of noise
(i.e., $X=0$), the matrix $M$, with low rank $k$, can be recovered exactly using convex relaxation under some conditions on the sampling rate $\delta$ and the singular vectors. These conditions were improved in \cite{candes2010power} and \cite{recht2011simpler}, and the approach was also extended to the case of noisy observed entries \cite{candes2010matrix}. The proposed algorithms involves solving a convex program, which can be computationally expensive. If the rank $k$ of the matrix is known, $M$ can be recovered using simpler spectral methods. For example, in \cite{keshavan2010matrix}, the authors show that in absence of noise, $M$ can be reconstructed asymptotically accurately using $O(\delta k mn \log n )$ operations under the conditions that the rank $k$ does not depend on $n $ and $m$, $\delta m = \omega(1)$ and $\delta n = \omega(1)$. Again these results can be adapted to the presence of noise \cite{keshavan2009matrix}. In this paper, we improve the spectral method used in \cite{keshavan2010matrix} and \cite{keshavan2009matrix}, so that it becomes memory-efficient, and so that it has performance guarantees even if the rank $k$ of $M$ scales with $m$ and $n$.

\medskip
\noindent
{\bf Streaming algorithms.} \cite{clarkson2009numerical} proposes an algorithm to provide a rank-$k$ approximation of a fully observed matrix $A$, using 1-pass on the columns of $A$. The algorithm uses a random $m\times \ell$ Rademacher matrix $S$, with an appropriate choice of $\ell$, and outputs a rank-$k$ matrix $\hat{A}^{(k)}$ constructed from $A^\dagger S$ and $AA^\dagger S$. When setting $\ell = O(k\varepsilon^{-1}\log(1/\eta))$ which requires $O(k\varepsilon^{-1}(m+n)\log(1/\eta))$ memory space, it is shown that with probability at least $1-\eta$,
\begin{equation} 
\|A- \hat{A}^{(k)} \|_F \le (1+\varepsilon)\|A- \bar{A}^{(k)} \|_F,\label{eq:con1}
\end{equation}
where $\bar{A}^{(k)}$ is the optimal rank-$k$ approximation of $A$. We could think of applying this algorithm to our problem. If the observed matrix $A$ is $A={\cal P}_{\Omega}(M+X)$, it would make sense to estimate $M$ by $\frac{1}{\delta}\hat{A}^{(k)}$ where $\hat{A}^{(k)}$ is the output of the algorithm in \cite{clarkson2009numerical} applied to $A$. Indeed, it is easy to check that $\|M-\frac{1}{\delta}\bar{A}^{(k)} \|_F^2 = o(mn)$ (i.e., the optimal rank-$k$ approximation of ${1\over \delta}A$ estimates $M$ asymptotically accurately). However, in general, $\frac{1}{\delta}\hat{A}^{(k)}$ is not asymptotically accurate:
\begin{eqnarray*}
\frac{\|M -\frac{1}{\delta}\hat{A}^{(k)} \|_F^2}{mn} &\ge& \frac{(\|A-
  \hat{A}^{(k)} \|_F - \|A- \bar{A}^{(k)}\|_F-\|\bar{A}^{(k)}- \delta
  M \|_F)^2}{\delta^2 mn}\cr
&= &\frac{(\varepsilon \|A- \bar{A}^{(k)}\|_F-\|\bar{A}^{(k)}- \delta
  M \|_F)^2}{\delta^2 mn}.
\end{eqnarray*}
Now, one can also easily check that $\|A- \bar{A}^{(k)} \|_F = \Theta (\sqrt{\delta mn})$ and $\|\delta M-\bar{A}^{(k)} \|_F = o(\delta \sqrt{mn})$, so that if we choose $\epsilon=\sqrt{\delta}$, we get $\frac{\|M -\frac{1}{\delta}\hat{A}^{(k)} \|_F^2}{mn}=\Omega(1)$. As a consequence, using the algorithm in \cite{clarkson2009numerical}, we cannot reconstruct $M$ asymptotically accurately using $O(k\sqrt{1/\delta}(m+n)\log(1/\eta))$ memory space. Recall that our algorithm reconstructs $M$ accurately with $O(k(m+n))$ memory space.  

We could also think of using sketching and streaming PCA algorithms to reconstruct $M$. When the columns arrive sequentially, these algorithms identify the left singular vectors in 1-pass on the matrix. We would then need a second pass on the data to estimate the right singular vectors, and complete the matrix. For example, \cite{liberty2013simple} proposes a sketching algorithm that updates the $\ell$ most frequent directions when a new column of $A$ is observed. This algorithm outputs a sketch $\hat{A}$ of $A$ and has the following performance guarantee: $\|A A^\dagger - \hat{A}\hat{A}^{\dagger} \|_2 \le \frac{2\|A\|_F^2}{\ell}$. It also uses $O(m\ell)$ memory space. Again if we apply the algorithm to our matrix completion problem, i.e., to the observed matrix $A={\cal P}_{\Omega}(M+X)$, where $M$ is of rank $k$, then $\|A \|_F^2 = \Theta(\delta m n)$ and $\sigma_k (AA^\dagger ) =\Theta( \delta^2 \sigma_k^2 (M)) = \Theta
(\frac{\delta^2 mn}{k})$. Hence to efficiently extract the top $k$ left singular vectors, we would need that $\frac{2\|A \|_F^2}{\ell} = o(\sigma_k (AA^\dagger ))$, which implies $\ell = \omega(k/\delta)$. Therefore, the required memory space would be $O(\frac{km}{\delta} + kn)$. Our algorithm is more efficient, and uses only 1-pass on the matrix. Note that the streaming PCA algorithm proposed in \cite{mitliagkas2013} does not apply to our problem (in \cite{mitliagkas2013}, the authors consider the spiked covariance model where a column is randomly generated in an i.i.d. every time).

 
\medskip
\noindent
{\bf Low complexity algorithms.} There have been recently an intense research effort to propose low-complexity algorithms for various linear algebra problems. Randomization has appeared as an efficient way to reduce the complexity of algorithms, see \cite{halko2011finding} for a survey. For example, \cite{sarlos2006improved} and \cite{clarkson2009numerical} devise algorithms for rank-$k$ approximation with guarantees (\ref{eq:con1}) and that use $O(\delta mn (k/\varepsilon + k \log k) + n\mbox{poly}(k/\varepsilon))$ operations. When the input matrix is sparse, \cite{clarkson2013low} leverages sparse embedding techniques, and reduces the required complexity to $O(\delta mn) + O((nk^2\varepsilon^{-4}+k^3\varepsilon^{-5})\cdot\mbox{polylog}(m+n))$ operations. But once again, as explained above, these results do not apply to our framework ((\ref{eq:con1}) is not enough to guarantee an asymptotically accurate matrix completion).



\begin{algorithm}[t!]
   \caption{Spectral PCA (SPCA)}
   \label{alg:pca}
\begin{algorithmic}
\STATE {\bfseries Input:} $A\in [0,1]^{m\times \ell}$, 
 $k$
\STATE $\hat{\delta} \leftarrow \frac{1}{m\ell}\sum_{(i,j)}1([A^{(B)}]_{ij}>0) $
\STATE {(Trimming)} $\tilde{A}\leftarrow$ erase rows of $A$ with more than $\max\{ 10, 10\hat{\delta} \ell \}$ non-zero entries 
\STATE $\Phi \leftarrow \tilde{A}^\dagger \tilde{A} - \mbox{diag}(\tilde{A}^\dagger \tilde{A})$
   \STATE $\hat{V}_{1:k} \leftarrow$ QR $(\Phi,k)$
   \STATE {\bfseries Output:} $\hat{V}_{1:k}$
\end{algorithmic}
\end{algorithm}
\begin{algorithm}[t!]
   \caption{QR Algorithm}
   \label{alg:pm}
\begin{algorithmic}
\STATE {\bfseries Input:} $\Phi$ (of size $\ell\times\ell$), $k$
\STATE {\bfseries Initialization:} $Q^{(0)} \leftarrow$ Randomly choose $k$
orthonormal vectors 
\FOR{$\tau=1$ {\bfseries to} $\lceil 10 \log (\ell) \rceil$}
\STATE $Q^{(\tau)} R^{(\tau)} \leftarrow$ QR decomposition of $\Phi Q^{(\tau-1)}$
\ENDFOR
\STATE {\bfseries Output:} $Q^{(\tau)}$
\end{algorithmic}
\end{algorithm}

\section{Extracting Right-Singular Vectors}

As mentioned in the introduction, the SMC algorithm deals with batches of arriving columns. Information from each batch will be extracted and aggregated as more columns arrive. In this section, we present an algorithm that will be used as a building block for extracting
information from a batch of columns. For concreteness, let assume that the size of a batch is $\ell$. In the SMC algorithm, $\ell$ will be chosen much smaller than $m$, so as to guarantee that the algorithm does not require large memory space. 

The algorithm presented in this section addresses the following problem. Let $M \in [0,1]^{m\times \ell}$ with singular value decomposition $M = U\Sigma
V^{\dagger}$. Given $0<k\leq \ell$ and $A=\mathcal{P}_{\Omega}(M+X)$, we wish to estimate the $k$ dominant right-singular vectors of $M$, $V_{1:k}$. At first, this might appear as a standard PCA task, but we are only interested in cases where $A$ is very sparse. Indeed $A$ only has a vanishing proportion $\delta$ of non-zero entries. Note that on average, we have $\delta \ell$ observed entries per row of $M+X$. Moreover, as this will become clear in the design of the SMC algorithm, we need to consider the case where $\delta \ell =o(1)$. In particular, there are many rows of $A$ with no observed entry. As a consequence, we do not get any information about the corresponding rows of $U$ in the singular value decomposition of $M$. Hence, we are here only interested in providing an estimate of the right-singular vectors $V$.

The algorithm to extract the dominant right-singular vectors, referred to as SPCA (Spectral Principal Component Analysis), is simple and its design relies on the following observation.
If we had access to the matrix $M$, then estimating the right-singular vectors of $M$ would be obvious. Indeed $M^{\dagger}M=V\Sigma^2 V^{\dagger}$, so that a standard
QR algorithm would output $V$. Now $A$ constitutes a subsampled noisy version of $M$ and we could try to apply this algorithm directly to $A$. From basic random matrix theory, we expect that the eigenvalues associated to the signal (i.e., the subsampled version of $M^{\dagger}M$) to be of the order of $\delta^2s_k^2(M)$. On the other hand, the eigenvalues associated with the noise (i.e., the subsampled version of  $X^{\dagger}X$) should be of the order $\delta \sqrt{m\ell}$. Thus, one could believe that the eigenvectors obtained by applying the QR algorithm to $A$ provide a good estimate of $V_{1:k}$ as soon as the ratio $\frac{\delta\sigma_k^2(M)}{\sqrt{m\ell}}$ is large enough. However, this is not quite true, because of the sparsity of the matrix $A$. To overcome this issue, we need to {\it regularize} the matrix $A$ before applying the QR algorithm. This is done in two steps:
\begin{itemize}
\item[(a)] {\em Trimming}: The rows of the subsampled matrix $A$ with too many non-zero entries are first removed. This trimming step is standard and avoids rows with too many
entries to perturb the spectral decomposition.
\item[(b)] {\em Removing diagonal entries:} Let $\tilde{A}$ denote the trimmed matrix. The diagonal entries of the covariance matrix $\tilde{A}^\dagger \tilde{A}$ are then removed: $\Phi = \tilde{A}^\dagger \tilde{A} -\mbox{diag}(\tilde{A}^\dagger \tilde{A})$. This step is needed because  the diagonal entries of $A^\dagger A$ scale as $\delta$, whereas its off-diagonal entries scale as $\delta^2$. Hence, when $\delta \to 0$, if the diagonal entries are not removed, they would be clearly dominant in the spectral decomposition. 
\end{itemize} 

In summary, the SPCA algorithm consists in applying the QR algorithm to the regularized version of $A$, i.e., to $\Phi$. Its pseudo-code is presented in Algorithm \ref{alg:pca}. The following theorem provides a performance analysis of SPCA, and is of independent interest.

\begin{theorem}
Let $\ell<m$, $\ell=o(1/\delta)$, and $M \in [0,1]^{m\times \ell}$ with singular value decomposition $M = U\Sigma
V^{\dagger}$, where $\Sigma = \mathrm{diag}(s_1(M),\dots, s_\ell(M))$ with
$s_1(M)\geq \dots\geq s_\ell(M)\geq 0$. Let $A= \mathcal{P}_{\Omega}(M+X)$.
Assume that there exists $k\leq \ell$ such that 
$s_k(M)=\omega(\sqrt{m})$, $\frac{s_k (M)}{s_{k+1}(M)} = \omega(1)$, and $\frac{\delta s_k^2(M)}{\sqrt{m\ell\log \ell}} =
\omega(1)$. Let $\hat{V}_{1:k}$ be the output of SPCA with input
$A$ and $k>0$. Then we have
$ \| (V_{1:k})^\dagger \cdot (\hat{V}_{1:k})_\bot \|_2 = o(1) $ with high
  probability. \label{pro:pca}
\end{theorem}

Note that the condition $\frac{\delta s_k^2(M)}{\sqrt{m\ell \log \ell}} =
\omega(1)$ in Theorem \ref{pro:pca} is similar to that suggested by the random matrix theory argument presented above. However we loose a $\log$ factor here because we use, in the proof, the Matrix Bernstein inequality (Theorem 6.1 of \cite{tropp2012user}). The condition $\frac{s_k (M)}{s_{k+1}(M)}\to \infty$ ensures a good separation in the spectrum of $M$ and is needed to ensure that the space spanned by $V_{k+1:\ell}$ is nearly orthogonal to the space
spanned by $V_{1:k}$ by Davis-Kahan $\sin \Theta$ Theorem (Theorem VII.3.2 in \cite{bhatia}). We conclude this section by analyzing the memory required by the SPCA algorithm, and its computational complexity.

\smallskip
\noindent{\bf Required memory.} SPCA needs to store $A$, $\Phi=\tilde{A}^\dagger \tilde{A} - \mbox{diag}(\tilde{A}^\dagger \tilde{A})$, and $\hat{V}$. The number of non-zero entries of $A$ is $O(\delta m\ell)$, and for each entry we need to store its id and its value. Hence for $A$, $O(\delta m\ell \log(m))$ memory is required. Similarly, the required memory for $\Phi$ is $O(\delta^2 m \ell^2 \log(\ell))$. Finally, storing $\hat{V}_{1:k}$ requires $O(\ell k)$ memory.
Overall the required memory is $O(\delta m\ell \log(m) + \ell k)$.

\smallskip
\noindent{\bf Computational complexity.} To run SPCA, we have to compute $\Phi$ and apply the QR algorithm to $\Phi$. The computation of $\Phi$ requires to perform $\frac{\ell(\ell-1)}{2}$ inner products of columns of $\tilde{A}$. Each inner product requires $O(\delta^2 m)$
floating-point operations, and thus the computational complexity to compute $\Phi$ is
$O(\delta^2m \ell^2)$. Now in the QR algorithm, we compute $\Phi Q_\tau$ and run the QR
decomposition $\log (\ell)$ times. The matrix product $\Phi Q_\tau$ requires $O(\delta^2m \ell^2 k )$ floating-point operations, while the QR decomposition requires
$O(\ell k^2)$ operations. Hence, the QR algorithm needs $O(\ell k (\delta^2 m\ell + k) \log (\ell) )$ operations. Overall, the computational complexity of SPCA is $O(\ell k (\delta^2 m\ell + k) \log (\ell) )$.

\section{Matrix completion with Streaming Input}

\begin{algorithm}[t!]
   \caption{Streaming Matrix Completion (SMC)}
   \label{alg:streaming1}
\begin{algorithmic}
\STATE {\bfseries Input:} $\{A_1,\dots,A_n \}$, 
$k$, $\ell$
\STATE 1. $A^{(B)} \leftarrow [A_1,\dots, A_{\ell}]$
\STATE 2. $\hat{\delta} \leftarrow \frac{1}{m\ell}\sum_{(i,j)}1([A^{(B)}]_{ij}>0) $ 
\STATE 3. $A^{(B_1)},A^{(B_2)},A^{(B_3)},A^{(B_4)} \leftarrow$ Split($A^{(B)},4,4,\hat{\delta}$)
\STATE 4. (PCA for the first block)$Q \leftarrow$ SPCA($A^{(B_1)},k$)
\STATE 5. (Trimming rows and columns)
\STATE $\quad\quad A^{(B_2)} \leftarrow$ make the rows having more than two
observed entries to zero rows
\STATE $\quad\quad A^{(B_2)} \leftarrow$ make the columns having more than $10m\hat{\delta}$
non-zero entries to zero columns
\STATE 6. (Reference Columns) $W \leftarrow A^{(B_2)}
Q$
\STATE 7. (Principle row vectors) $\hat{V}^{1:\ell} \leftarrow
(A^{(B_3)})^{\dagger} W $  
\STATE 8. (Principle column vectors) $\hat{I} \leftarrow A^{(B_4)} 
\hat{V}^{1:\ell}$
\STATE {\em Remove $A^{(B)}$, $A^{(B_1)}$, $A^{(B_2)}$, $A^{(B_4)}$,
  and $Q$ from the memory space}
\FOR{$t =\ell+1$ {\bfseries to} $n$}
\STATE 9. $A^{(1)}_t,A^{(2)}_t \leftarrow$ Split($A_t,2,4,\hat{\delta}$)
\STATE 10. (Principle row vectors) $\hat{V}^{t} \leftarrow
(A^{(1)}_t)^{\dagger} W$
\STATE 11. (Principle column vectors) $\hat{I} \leftarrow \hat{I}  +
A^{(2)}_t \hat{V}^{t}$
\STATE {\em Remove $A_t$ and $A'_t$ from the memory space}
\ENDFOR
\STATE 12. $\hat{R}\leftarrow$ find $\hat{R}$ using the Gram-Schmidt
process such that $\hat{V}\hat{R}$ is an orthonormal matrix.
\STATE 13. $\hat{U}\leftarrow  \frac{4}{\hat{\delta}} \hat{I}\hat{R}\hat{R}^{\dagger}$
\STATE {\bfseries Matrix completion:} $|\hat{U} \hat{V}^\dagger |^1_0$
\end{algorithmic}
\end{algorithm}

\begin{algorithm}[t!]
   \caption{Split}
   \label{alg:split}
\begin{algorithmic}
\STATE {\bfseries Input:} $A$,$a$,$b$, $\delta$
\STATE {\bfseries Initial:} $A^{(1)},\dots,A^{(a)} \leftarrow$ zero
matrices having the same size as $A$
\FOR {every $[A]_{uv}$}
\STATE $\gamma \leftarrow$ $s \subset \{1,\dots,b\}$ which is randomly selected over all subsets of $\{1,\dots,b\}$
with probability
$\frac{1}{\delta}\left(\frac{\delta}{b}\right)^{|s|}\left(1-\frac{\delta}{b}\right)^{b-|s|}$
if $s$ is not the empty set and with probability
$1-\frac{1}{\delta}(1-(1-\frac{\delta}{b})^b)$ if $s$ is the empty set
\FOR {$i \in \gamma$}
\STATE $[A^{(i)}]_{uv} \leftarrow [A]_{uv}$ 
\ENDFOR
\ENDFOR
\STATE {\bfseries Output:} $A^{(1)},\dots,A^{(a)}$
\end{algorithmic}
\end{algorithm}

In this section, we present our main algorithm, SMC, that reconstructs a matrix $M\in [0,1]^{m\times n}$ from a few noisy observations on its entries, i.e., from $A=\mathcal{P}_\Omega(M+X)$. The pseudo-code of SMC is presented in Algorithm~\ref{alg:streaming1}. SMC consists in three main steps: {\em Step 1)} Generate reference columns denoted by $W$, 
{\em Step 2)} Find principle row vectors $\hat{V}$ using $W$, and {\em
  Step 3)} Find $\hat{U}$ such that $\hat{U}\cdot \hat{V}^{\dagger}
\approx M$. In what follows, we explain each of these steps in details and
show for each step which conditions of Assumption 1 are needed. All proofs are presented in Appendix. The singular value decomposition of $M$ is $M=U\Sigma V^{\dagger}$.

\subsection{Step 1: Finding reference columns $W$}\label{sec:step1}

We now explain the first step of the algorithm leading to a $m\times
k$ matrix $W$ containing {\it reference} columns. This step corresponds to lines 1 to 6 in the pseudo-code. 

Let $A^{(B)}=A_{1:\ell}$ be the batch of the $\ell$ first arriving columns of
$A$. Note in particular that we have:
$$
A^{(B)} = \mathcal{P}_{\Omega}(M^{(B)} + X_{1:\ell}) \quad \mbox{
  with, } \quad M^{(B)} = M_{1:\ell}=U\Sigma \left( V^{1:\ell}\right)^\dagger.
$$
In line $2$, we compute $\hat{\delta}$, an estimate of the sampling rate $\delta$. In line $3$, we construct $4$ undersampled copies of $A^{(B)}$. For $i\in \{1,2,3,4\}$, the different $A^{(B_i)}$'s are independent given $M+X$ and have the same distribution as $A^{(B)}$, except that the parameter $\delta$ is now replaced by $\delta/4$. 

The first non-trivial operation is presented in line 4 where we apply the
algorithm SPCA described in previous section to the matrix
$A^{(B_1)}$. In order to apply our Theorem \ref{pro:pca}, we need to
have:
\begin{eqnarray}
\label{eq:condMB}\frac{s_k (M^{(B)})}{s_{k+1}
  (M^{(B)})} = \omega(1)\quad\mbox{and}\quad
\frac{\delta s_k^2(M^{(B)})}{\sqrt{m\ell \log
  \ell}} = \omega(1).
\end{eqnarray}
Note that there is a slight abuse of notation as the distribution of
$A^{(B_1)}$ is the same as the one of $A^{(B)}$ if we change $\delta$
to $\delta/4$ but a constant factor $4$ is clearly irrelevant here.
Our first task is to translate the conditions (\ref{eq:condMB}) on the
original matrix $M$.
To this aim, we state the following lemma:
\begin{lemma} Let $M=U\Sigma V^\dagger$ be a $m\times n$ matrix and
  $\ell\leq n$. Denote by $M^{(B)}=M_{1:\ell}$.
If $s_k^2 (M) = \omega(\frac{mn \log m}{\ell})$ and $\frac{s_k(M)}{s_{k+1}(M)}=\omega(1)$, then with high probability, 
$$s_k(U_{1:k}U_{1:k}^\dagger M^{(B)}) \ge
\sqrt{\frac{\ell}{2n}}s_k(M) \quad\mbox{and}\quad
\frac{s_k(U_{1:k}U_{1:k}^\dagger M^{(B)})}{s_1((I-U_{1:k}U_{1:k}^\dagger) M^{(B)})} = \omega(1).$$
\label{lem:cor}
\end{lemma}
Its proof is given in Appendix~\ref{sec:proof-lemma3} and follows from the matrix Chernoff bound (Theorem 2.2 of
\cite{tropp2011improved}). 

Note that $U_{1:k}U_{1:k}^\dagger$ is the orthogonal projection on the
span of $U_{1:k}$. As a result, we have $s_k(M^{(B)})\geq
s_k(U_{1:k}U_{1:k}^\dagger M^{(B)})$ by a simple application of the
Courant-Fischer variational formulas for singular values. In
particular, as soon as $\frac{\ell \delta^2 s_k^4(M)}{mn^2\log
  \ell}\to \infty$, we see that the second condition in
(\ref{eq:condMB}) is satisfied.
To get the first condition in (\ref{eq:condMB}), we write:
$$
M^{(B)}= U_{1:k}U_{1:k}^\dagger M^{(B)}+ (I-U_{1:k}U_{1:k}^\dagger) M^{(B)},
$$
note that the first matrix is of rank $k$ and we can use Lidskii's 
inequality $s_{k+1}(A+B) \leq s_k(A)+s_1(B)$ to get:
$$
s_{k+1}(M^{(B)}) \leq s_1((I-U_{1:k}U_{1:k}^\dagger) M^{(B)}).
$$
Hence we have
$$
\frac{s_k(M^{(B)})}{s_{k+1}(M^{(B)})}\geq\frac{s_k(U_{1:k}U_{1:k}^\dagger M^{(B)})}{s_1((I-U_{1:k}U_{1:k}^\dagger) M^{(B)})},$$
and the first condition in (\ref{eq:condMB}) follows from the second
statement in Lemma \ref{lem:cor} as soon as its conditions are
satisfied. Combined with Lemma \ref{lem:cor}, Theorem \ref{pro:pca}
allows us to get the properties of $Q$ computed in line 4 of the Algorithm SMC:
\begin{corollary}
Assume that there exists $k$ and $\ell$ such
that $\frac{s_k (M)}{s_{k+1} (M)} = \omega (1)$,
$\frac{\delta^2 \ell s_k^4(M)}{ m n^2 \log \ell} = \omega(1)$,
and $s_k^2 (M) = \omega(\frac{mn \log m}{\ell})$.
Let $\bar{V}^{1:\ell}$ be an orthonormal basis of the
  linear span of $V^{1:\ell}_{1:k}$. 
Then we have $\|(\bar{V}^{1:\ell})^\dagger \cdot Q_\bot\| =
o(1)$ with high probability, where $Q$ is the $\ell \times k$ matrix obtained in line 4 of the Algorithm SMC.
\label{cor:1}
\end{corollary}

Once we have $Q$, we compute what we call the reference columns as follows:
$$W =  A^{(B_2)} \cdot Q.$$ 
Note that $W$ will be kept in memory during the whole algorithm. It is
relatively easy to see that the linear span of the columns of $W$ is a
noisy version of the linear span of $U_{1:k}$.
Indeed, note that $\mathbb{E}[A^{(B_2)}] = \frac{\delta}{4} M^{(B)}$,
moreover we have $M^{(B)}=U\Sigma (V^{1:\ell})^{\dagger}\approx
U_{1:k}\Sigma_{[k]} (V^{1:\ell}_{1:k})^{\dagger}$ thanks to Lemma
\ref{lem:cor}. Hence we have
$$
W = A^{(B_2)} Q \approx \frac{\delta}{4}U_{1:k}\Sigma_{[k]} (V^{1:\ell}_{1:k})^{\dagger}Q.
$$
By Corollary \ref{cor:1}, the span of the columns of $Q$ is
approximately the span of the column of $V^{1:\ell}_{1:k}$ so that the
singular values associated to the linear span of $U_{1:k}$ are
$\Omega(\delta s_k(M^{(B)}))=\Omega(\delta \sqrt{\ell/n}s_k(M))$ by
Lemma \ref{lem:cor}.
This value has to be compared to the noise level. For the same reason
as in Section 3, we first trim the matrix $A^{(B_2)}$
(note that the first trimming phase in line 5 is made to ease the
technical proof).
After the trimming process, the singular values of $(A^{(B_2)}-\mathbb{E}[A^{(B_2)}] )\cdot
 Q$ are bounded by $O(\sqrt{\delta m\ell})$. Unfortunately, in our setting this can be
 much larger than $\delta \sqrt{\ell/n}s_k(M)$. However, the hidden
 signal in $W$ is in the span of the columns of $U_{1:k}$ and all the
 columns that arrive belong (approximately) to this span. In the
 sequel, we use this fact in order to amplify the signal in $W$ when estimating $V$ and then $U$.

\subsection{Step 2: Finding principle row vectors $\hat{V}$}\label{sec:step2}

In this section, we explain how we recover $V_{1:k}$ or at least $k$ vectors having
the same linear span as $V_{1:k}$. 

Let $A^{(1)} = [A^{(B_3)},A^{(1)}_{\ell+1},\dots,A^{(1)}_n]$. 
Note that thanks to the splitting procedure in line 9, the columns of
$A^{(1)}$ are i.i.d. with sampling rate $\delta/4$.
In the SMC algorithm, we simply get an estimate of $V$ as follows: $\hat{V} =
(A^{(1)})^{\dagger} W$. 
The linear span of the columns of $\hat{V}$ becomes
very close to the linear span of the columns of $V_{1:k}$ when
$$\frac{s_k( V_{1:k} V_{1:k}^\dagger \hat{V} )}{ s_1((I- V_{1:k} V_{1:k}^\dagger) \hat{V})} = \omega(1) .$$
This can be seen as in Section \ref{sec:step1} since $V_{1:k}
V_{1:k}^\dagger$ is simply the orthogonal projection on $V_{1:k}$.

The above condition holds for the following reasons:
\begin{itemize}
\item The signal is amplified (Lemma~\ref{lem:pf:thm3-1} in
  Appendix). Since $\mathbb{E}[A^{(1)}] = \frac{\delta}{4} M$, we see
  that
\begin{eqnarray*}
\hat{V}=(A^{(1)})^\dagger W &\approx&  \frac{\delta^2}{16}V\Sigma U^\dagger
U_{1:k}\Sigma_{[k]} (V^{1:\ell}_{1:k})^{\dagger}Q\\
&\approx&\frac{\delta^2}{16}V_{1:k}\Sigma_{[k]}^2 (V^{1:\ell}_{1:k})^{\dagger}Q.
\end{eqnarray*}
Roughly, the signal which was $\Omega(\delta \sqrt{\ell/n}s_k(M))$ is
now multiplied by $\delta s_k(M)$ and we get:
$$s_k( V_{1:k} V_{1:k}^\dagger \hat{V} )= \Omega(V_{1:k}
V_{1:k}^\dagger (\mathbb{E}[A^{(1)}])^\dagger \mathbb{E}[A^{(B_2)}]Q ) =  \Omega(\delta^2 s^2_k(M)  \sqrt{\frac{\ell}{n}}).$$
\item The noise is cancelled (Lemma~\ref{lem:thm3-2} in Appendix). Since the two noise matrices $A^{(B_2)}-\mathbb{E}[A^{(B_2)}]$
  and  $A^{(1)}-\mathbb{E}[A^{(1)}]$ are independent, the noise
  directions are not amplified as much as the signals. We can bound
  the noise as  follows:
$$s_1((I- V_{1:k} V_{1:k}^\dagger) \hat{V} ) =
o(\delta^2 s^2_k(M)  \sqrt{\frac{\ell}{n}}) .$$
\end{itemize}
Putting things togetehr, we obtain the following result:
\begin{theorem}
Assume that there exists $k$ and $\ell$ such
that $s_k^2 (M) = \omega(\frac{mn \log m}{\ell})$, $\frac{s_k (M)}{s_{k+1} (M)} = \omega (1)$, and $\frac{\delta^2 \ell s_k^4
  (M)}{ m n^2(k+ \log \ell)} = \omega(1)$. Then we have $\|V_{1:k}^{\dagger} (\hat{V}_{1:k})_{\bot} \| = o(1)$ with high probability.\label{thm:streamingV}
\end{theorem}

\subsection{Step 3: Finding principle column vectors $\hat{U}$}
In the previous step, we identified a $n \times
k$ matrix $\hat{V}$ estimating the principle row vectors of $M$. From this estimate, we now extract the matrix $\hat{U}$ such that
$\|\hat{U}\hat{V}^\dagger - M\|_F = o(mn)$.

Let $A^{(2)} = [A^{(B_4)},A^{(2)}_{\ell+1},\dots,A^{(2)}_n]$.
For simplicity, suppose that the
linear span of the rows of $\hat{V}^\dagger$ is exactly the same as the linear span
of the rows of $M$. From $\hat{V}$, we can generate a $k\times k$
matrix $\hat{R}$ using the Gram-Schmidt process so that
$ \hat{V} \hat{R} $ becomes an orthogonal matrix. Since
$\hat{V} \hat{R}$ is an orthonormal basis of the linear span of the
rows of $M$, we have
$$M = \frac{4}{\delta} \mathbb{E}[A^{(2)}]\hat{V}\hat{R}(\hat{V}\hat{R})^\dagger =
(\frac{4}{\delta} \mathbb{E}[A^{(2)}]\hat{V}\hat{R}\hat{R}^{\dagger})
\cdot \hat{V}^\dagger = \bar{U}\hat{V}^\dagger,$$ 
where $\bar{U} = \frac{4}{\delta}\mathbb{E}[A^{(2)}]\hat{V}\hat{R}\hat{R}^{\dagger}$.

From the above observation, we propose to compute $\hat{U}$ as follows:
\begin{eqnarray*}
\hat{U}&=&\frac{4}{\delta}\hat{I}\hat{R} \hat{R}^\dagger ~ = ~
\frac{4}{\delta} A^{(2)} \hat{V}\hat{R} \hat{R}^\dagger \cr
&= &\frac{4}{\delta} \mathbb{E}[A^{(2)}] \hat{V}\hat{R} \hat{R}^\dagger +\frac{4}{\delta} (A-\mathbb{E}[A^{(2)}]) \hat{V}\hat{R} \hat{R}^\dagger .
\end{eqnarray*}
Then, we need to prove that 
the row space of $A^{(2)}-\mathbb{E}[A^{(2)}]$ is almost orthogonal to
$V$, to get $\hat{U} = \bar{U}+(A-\mathbb{E}[A]) \hat{V}\hat{R}
\hat{R}^\dagger = (1+o(1))\bar{U}$. This is true only if $n$ is large
enough, indeed $n=\omega(k/\delta)$ (see Appendix).

We are now ready to analyze the performance of the SMC algorithm. 
We first need to check that Assumption 1 implies the technical
conditions required in our previous results.
When $M$ has $k$ dominant singular values such that
$\frac{s_k(M)}{s_{k+1}(M)} = \omega(1)$ and
$\sum_{i>k} s_{i}^2(M) = o(mn)$, then, $s^2_k(M) =
\Omega(\frac{mn}{k})$.
To see this, assume this is not the case so that there exists $k' < k$ such that
$s_{k'} = \Omega(\frac{mn}{k})$ and $s^2_{k'+1}(M) =
o(\frac{mn}{k})$. But then $\frac{s_{k'}(M)}{s_{k'+1}(M)} = \omega(1)$
and $\sum_{i>k'}s^2_{i} = o(mn)$ which contradicts
the minimality of $k$. Therefore, the conditions $s_k^2 (M) = \omega(\frac{mn \log m}{\ell})$ and $\frac{\delta^2 \ell s_k^4
  (M)}{ m n^2(k+ \log \ell)} = \omega(1)$ become $\ell = \omega(k \log m)$ and $\frac{\delta^2 m \ell}{
  k^2 (k+ \log \ell)} = \omega(1)$, which are satisfied by Assumption
1 when $\ell =O( m)$ and
$\ell = \Omega(\frac{k}{\delta \log m})$.
Hence we obtain the following result:
\begin{theorem}
Assume that Assumption 1 is satisfied with $\ell =
\Omega(\frac{k}{\delta \log m})$ and $\ell=O(m) $. Then with high probability, the SMC algorithm provides an asymptotically accurate estimate of $M$:
$$\frac{\|M-[\hat{U}\hat{V}^\dagger]^1_0\|_F}{mn} = o(1). $$\label{thm:stmcomp}
\end{theorem}

\subsection{Required Memory}

Next we analyze the memory required by the SMC algorithm. \\
{\it From line 1 to 8 in the pseudo-code.} We need to store $A^{(B)}$, $A^{(B_1)}$, $A^{(B_2)}$, $A^{(B_3)}$, and $A^{(B_4)}$. Since these matrices are sparse with sampling rate $\delta$ or $\delta/4$, we need to store only $O(\delta m \ell)$ of their elements and $O(\delta m \ell \log m)$bits to store the id of the non-zero entries. From the previous section, we know that the {\em SPCA} algorithm requires $O(\delta m \ell \log m + k \ell)$ memory to find $Q$. Finally we need to store $\hat{V}$ and $\hat{I}$.
Thus, when $\ell =\frac{k}{\delta \log m}$, this first part of the algorithm requires $O(km+kn)$.\\
{\it From line 9 to 11.} Here we treat the remaining columns. Note that before doing that, $A^{(B)}$, $A^{(B_1)}$, $A^{(B_2)}$, $A^{(B_3)}$,
  and $Q$ are removed from the memory. Using this memory, for the $t$-th arriving column, we can store it, compute $\hat{V}^{t}$ and $\hat{I}$, and remove the column to save memory. Therefore, we do not need additional memory to treat the remaining columns.\\
{\it Lines 12 and 13.} From $\hat{I}$ and $\hat{V}$, we compute $\hat{U}$. To this aim, the memory required is $O(km+kn)$.

In summary, we have:
\begin{theorem}
When $\ell = \frac{k}{\delta \log (m)}$,
the memory required to run the SMC algorithm is $O(km+kn)$.\label{thm:memory}
\end{theorem}

\subsection{Computational Complexity}
The computational complexity of the {\em SMC} (Algorithm~\ref{alg:streaming1}) depends on the number of non-zero elements of $A$ and $\ell$. More precisely:\\
{\it From line 1 to 8.} From the previous section, the {\em SPCA} algorithms requires $O(\ell
k(\delta^2 m\ell + k) \log (\ell) )$ floating-point operations to compute
$Q$. The computations of $W$, $\hat{V}$, and $\hat{I}$ are
just inner products, and require $O(\ell k(\delta^2 m\ell + k) \log (\ell) )$ operations. \\
{\it From line 9 to 11.} To compute
  $\hat{V}^{t}$ and $\hat{I}$ when the $t$-th column arrives, we need $O(k m \delta)$ operations. Since there are $n-\ell$ remaining
  columns, the total number of operations is $O(k m n \delta)$.\\
{\it Lines 12 and 13} $\hat{R}$ is computed from $\hat{V}$ using the
  Gram-Schmidt process which requires $O(k^2m)$ operations. We then compute $\hat{I}\hat{R}\hat{R}^{\dagger}$ using $O(k^2m)$ operations .

When $\ell = \frac{k}{\delta \log (m)}$ and $k^2 = O(\delta n)$, the number of operations to treat the first $\ell$ columns is
\begin{align*}O(\ell k(\delta^2 m\ell + k) \log (\ell) ) &=O(k\delta^2 m\ell^2 \log
(\ell)) +O(\ell k^2 \log (\ell) ) \cr &= O(k^3 m \frac{\log \ell}{\log^2
  m}) + O(\delta mn)= O(k mn \delta) .
\end{align*}
Since the remaining part of the algorithm requires $O(\delta k m n)$ operations as well, we conclude:
Theorem~\ref{thm:complex}.
\begin{theorem}
Assume that Assumption 1 is satisfied with $\ell = \frac{k}{\delta
  \log (m)}$. Then, the computational complexity of the SMC algorithm is $O(\delta k m n)$. \label{thm:complex}
\end{theorem}

\section{Conclusion}
This paper investigated the streaming memory-limited matrix completion problem
when the observed entries are noisy versions of a small random
fraction of the original entries. We proposed a streaming algorithm
which produces an estimate of the original matrix with
a vanishing mean square error, uses memory space scaling linearly with the ambient dimension of the matrix, i.e. the memory required to store the output alone, and spends
computations as much as the number of non-zero entries of the input matrix. Our algorithm is relatively simple, and in particular, it does exploit elaborated techniques (such as sparse embedding techniques) recently developed to reduce the memory requirement and complexity of algorithms addressing various problems in linear algebra.

\newpage
\bibliography{reference}

\newpage
\appendix

\section{Appendix}

\subsection{Proof of Theorem~\ref{pro:pca}}
We can split $\Phi$ as follows:
\begin{align*}
\Phi =& \delta^2 V_{1:k}V_{1:k}^\dagger M^\dagger M + \Phi - \delta^2 V_{1:k}V_{1:k}^\dagger M^\dagger M.
\end{align*}
The power method can find $\hat{V}$ such that $\|\hat{V}^{\dagger} (V_{1:k})_{\bot}
\|_2 = o(1)$ when
$\frac{\delta^2  s_k (M^{\dagger}M)}{\|\Phi - \delta^2 V_{1:k}V_{1:k}^\dagger M^\dagger M \|_2} =
\omega(1)$ which is shown in Lemma 11 of \cite{yun2014nips}. Since 
\begin{eqnarray}
\|\Phi - \delta^2 V_{1:k}V_{1:k}^\dagger M^\dagger M\|_2 &\le& \|\mathbb{E}[\Phi]- \delta^2 V_{1:k}V_{1:k}^\dagger M^\dagger M \|_2 + \|\Phi - \mathbb{E}[\Phi]
\|_2 \cr
&\le& \|\delta^2   \mbox{diag}(M^{\dagger}M)\|_2 + \| \delta^2 (I-V_{1:k}V_{1:k}^\dagger) M^\dagger M\| + \|\Phi - \mathbb{E}[\Phi]
\|_2 \cr &\le& \delta^2 m + \delta^2 s^{2}_{k+1}( M) + \|\Phi - \mathbb{E}[\Phi]
\|_2,\label{eq:noise}\end{eqnarray}
in the remaining part, we transform $\Phi-\mathbb{E}[\Phi]$ as a
sum of random matrices, and then using Matrix Bernstein
inequality we get an upper bound for $\|\Phi - \mathbb{E}[\Phi]
\|_2$ to conclude this proof.

Recall that $A^i$ is the $i$-th low of $A$ and 
\begin{align*}
\Phi-\mathbb{E}[\Phi] & = \sum_{i=1}^{m} \left((A^i)^\dagger A^i
  - \mbox{diag}((A^i)^\dagger A^i) -
  \mathbb{E}[(A^i)^\dagger A^i - \mbox{diag}((A^i)^\dagger A^i)]  \right) .
\end{align*}
Let $X^{(i)} = (A^i)^\dagger A^i - \mbox{diag}((A^i)^\dagger A^i) -
  \mathbb{E}[(A^i)^\dagger A^i - \mbox{diag}((A^i)^\dagger A^i)].$ 
Then $X^{(i)}$ is a self-adjoint $\ell\times \ell$ matrix and $\mathbb{E}[X^{(i)}] =
0$.

The Matrix Bernstein inequality (Theorem 6.1 \cite{tropp2012user}) is
a matrix concentration inequality for the sum of zero mean random matrices.
\begin{proposition}[Matrix Bernstein] Consider a finite independent
  random matrix set $\{ X^{(i)} \}_{1\le i \le m}$, where every $X^{(i)}$ is
  self-adjoint with dimension $n$, $\mathbb{E}[X^{(i)}] = 0$, and
  $\|X^{(i)} \|_2 \le R$ almost surely. Let $\rho^2 = \| \sum_{i=1}^{m}
  \mathbb{E}[X^{(i)}  X^{(i)}] \|_2$. Then,
$$\mathbb{P} \{ \| \sum_{i=1}^{m} X^{(i)} \|_2 \ge x \} \le n
\exp\left(\frac{- x^2/2}{\rho^2 + Rx/3} \right).$$\label{pro:be}
\end{proposition} 
In order to use the Matrix Bernstein inequality,
we have to find upper bounds for $\|X^{(i)}\|_2$ and $\rho^2 $. Since $A^i$ are independently sampled with
  probability $\delta$, $[X^{(i)}]_{uv}$ has a some constant value if both $u$
  and $v$ are sampled in $A_i$ and $O(\delta^2)$ otherwise. Using these, the
  following lemmas bound $\|X^{(i)}\|_2$ and $\rho^2 $. 
\begin{lemma}
When $n=\omega(1)$, for $1\le i \le m$, there
exists a constant $C_1$ such that
$$\|X^{(i)}\|_2  \le C_1 \max\{ 1, \delta \ell \}.$$\label{lem:b1}
\end{lemma}
\noindent{\bf \em Proof:} Since the number of non-zero entries of
$A^i$ is bounded by $\max\{10, 10\delta \ell \}$, we can easily compute $r_u = \sum_{v\neq
  u} |[X^{(i)}]_{uv}| \le \max\{10, 10\delta \ell \} + \delta \ell$
for all $1\le i \le m$ and $1\le u \le \ell$. By the Gershgorin circle
theorem, therefore, for all $i$
$$\|X^{(i)}\|_2 \le \max\{10, 10\delta \ell \} + \delta \ell.$$
\QED
\begin{lemma}There exists a constant $C_2$ such that
$$\| \sum_{i=1}^{m} \mathbb{E}[X^{(i)}X^{(i)}] \|_2  \le  C_2m \max\{
\delta^2 \ell
, \delta^3 \ell^2 \}.$$ \label{lem:b2}
\end{lemma}
\noindent{\bf \em Proof:} Since the number of non-zero entries of
$A^i$ is bounded by $\max\{10, 10\delta \ell \}$, every
$|\mathbb{E}[X^{(i)}X^{(i)}]_{uv}|=O(\delta^2 (1+\delta\ell)) $ when
$u\neq v$ and every $|\mathbb{E}[X^{(i)}X^{(i)}]_{uu}|=O(\delta^2 \ell
(1+\delta \ell))$. By the Gershgorin circle theorem, therefore
$$\| \sum_{i=1}^{m} \mathbb{E}[X^{(i)}X^{(i)}] \|_2 = O(\delta^2m \ell
(1+ \delta \ell)).$$
\QED

Let $C = 16  \max\{C_1 , C_2 \}$. From Lemma~\ref{lem:b1} and \ref{lem:b2} and Proposition~\ref{pro:be},
\begin{equation}
\mathbb{P} \left\{\|\Phi - \mathbb{E}[\Phi]
\|_2 \ge \sqrt{C \log (n)  
  \max\{1, \delta^2m\ell, \delta^3 m \ell^2 \}} \right\} \le \frac{1}{\ell^2} .\label{eq:xbnd}
\end{equation}

\vspace{0.3cm}
\noindent{\bf \em Proof of Theorem~~\ref{pro:pca}:} This proof starts
with 
$$
\Phi = \delta^2 V_{1:k}V_{1:k}^\dagger M^\dagger M + \Phi - \delta^2
V_{1:k}V_{1:k}^\dagger M^\dagger M = \delta^2 V_{1:k}V_{1:k}^\dagger M^\dagger M +Y,
$$
where $Y=\Phi - \delta^2 V_{1:k}V_{1:k}^\dagger M^\dagger M $. From \eqref{eq:noise} and \eqref{eq:xbnd}
\begin{align*}\|Y
\|_2 &\le \delta^2 m + \delta^2 s^{2}_{k+1}( M) + \sqrt{C \log (\ell)   \max\{1,
  \delta^2m \ell, \delta^3 m \ell^2 \}} \cr
&= o(\delta^2 s^2_{k}(M)) + \sqrt{C \log (\ell)   \max\{1,
  \delta^2m \ell \}},\end{align*}
where the last equality stems from $s^2_{k}(M) = \omega(m)$ and
$\frac{s_k (M)}{s_{k+1}(M)} = \omega(1)$ the conditions of this theorem. Since the
condition $\frac{\delta^2 s^4_k(M)}{m\ell \log \ell } = \omega(1)$ implies
$\delta^2 m \ell = \omega(k^2 \log \ell)$ and
$\frac{\delta^2 s^2_k(M)}{\sqrt{C \log (\ell)   \max\{1,
  \delta^2m \ell \}}}=\frac{\delta^2 s^2_k(M)}{\sqrt{C \log (\ell)  
  \delta^2m \ell}} = \omega(1)$, we can deduce 
$\frac{s_k(\delta^2 V_{1:k}V_{1:k}^\dagger M^\dagger M)
}{\|Y\|_2} = \omega(1)$. Therefore, $\|(V_{1:k})^\dagger (\hat{V})_\bot \| =
o(1)$ from Lemma 11 of  \cite{yun2014nips}.

\subsection{Proof of Lemma~\ref{lem:cor}} \label{sec:proof-lemma3}

Let $F = U_{1:k}^\dagger M^{(B)}$ and $G=
(I-U_{1:k}U_{1:k}^\dagger) M^{(B)}$. We find a lower bound for
$s_k(F)$ and an upper bound $s_1(G)$ using the matrix
Chernoff bound (Theorem 2.2 in \cite{tropp2011improved}).
\begin{proposition}[Matrix Chernoff]
Let $\mathcal{X}$ be a finite set of positive-semidefinite matrices
with dimension $d$ and satisfy $\max_{X\in \mathcal{X}} s_1 (X)
\le \alpha$. Let
$$\beta_{\min} = \frac{\ell}{|\mathcal{X}|} s_d (\sum_{X \in
  \mathcal{X}}X) \quad\mbox{and}\quad \beta_{\max} = \frac{\ell}{|\mathcal{X}|} s_1 (\sum_{X\in \mathcal{X}}X).$$
When $\{X^{(1)},\dots,X^{(\ell)}\}$ are sampled uniformly
at random from $\mathcal{X}$ without replacement,
\begin{align*}
\mathbb{P}\left\{ s_{1}(\sum_{i=1}^{\ell}X^{(i)}) \ge
  (1+\varepsilon) \beta_{\max} \right\} &\le d
  \left(\frac{e^{\varepsilon}}{(1+\varepsilon)^{1+\varepsilon}}
  \right)^{\beta_{\max}/\alpha}\quad\mbox{for}~\varepsilon \ge 0
                                          \quad\mbox{and}\cr
\mathbb{P}\left\{ s_{d}(\sum_{i=1}^{\ell}X^{(i)}) \le
  (1-\varepsilon) \beta_{\min} \right\} &\le d
  \left(\frac{e^{-\varepsilon}}{(1-\varepsilon)^{1-\varepsilon}}
  \right)^{\beta_{\min}/\alpha}\quad\mbox{for}~\varepsilon\in [0,1).
\end{align*}\label{pro:matrix}
\end{proposition}

\noindent{\em i) $s_k(F)$:} $FF^\dagger$ is the sum of $\ell$ matrices
which are sampled uniformly at random from 
$\mathcal{X} = \{ U_{1:k}^\dagger M_1(U_{1:k}^\dagger
M_1)^\dagger,\dots,U_{1:k}^\dagger M_n(U_{1:k}^\dagger
M_n)^\dagger \}$ without replacement where the matrix dimension is
$k$. We can obtain the other parameters to compute the matrix Chernoff as follows:
$\alpha = m$ since $\|M_i\|^2 \le m$ for all $1\le i \le n$ and
$\beta_{\min} = \frac{\ell}{n}s^2_k(M)$. From Proposition~\ref{pro:matrix},
$$\mathbb{P}\left\{ s_{k}(FF^\dagger) \le
  (1-\varepsilon) \frac{\ell}{n}s^2_k(M) \right\} \le k
  \left(\frac{e^{-\varepsilon}}{(1-\varepsilon)^{1-\varepsilon}}
  \right)^{\frac{\ell}{mn}s^2_k(M)}\quad\mbox{for}~\varepsilon\in [0,1).$$
Therefore, when $s^2_k (M) = \omega (\frac{mn\log m}{ \ell })$, 
$$\mathbb{P}\left\{ s^2_{k}(F) \le
  \frac{\ell}{2n}s^2_k(M) \right\} \le \frac{1}{m}.$$

\smallskip
\noindent{\em ii) $s_1(G)$:} $GG^\dagger$ is the sum of matrices
sampled uniformly at random without replacement from
$$\mathcal{X} = \{ (I-U_{1:k}U_{1:k}^\dagger) M_1((I-U_{1:k}U_{1:k}^\dagger)
M_1)^\dagger,\dots,(I-U_{1:k}U_{1:k}^\dagger) M_n((I-U_{1:k}U_{1:k}^\dagger)
M_n)^\dagger \}.$$ Here, the dimension is $m$, $\alpha = m$ and $\beta_{\max} = \frac{\ell}{n}s^2_{k+1}(M)$.
From Proposition~\ref{pro:matrix},
$$\mathbb{P}\left\{ s_{1}(GG^\dagger) \ge
  (1+\varepsilon) \frac{\ell}{n}s^2_{k+1}(M) \right\} \le m
  \left(\frac{e^{\varepsilon}}{(1+\varepsilon)^{1+\varepsilon}}
  \right)^{\frac{\ell}{mn}s^2_{k+1}(M)}\quad\mbox{for}~\varepsilon
  \ge 0.$$
When we set $\varepsilon^{\star} = \max\{2, \frac{2 m n \log m}{\ell
  s_{k+1}^2(M)} \}$, $\mathbb{P}\left\{ s_{1}(GG^\dagger) \ge
  (1+\varepsilon^{\star}) \frac{\ell}{n}s^2_{k+1}(M) \right\}
\le \frac{1}{m}$ and $(1+\varepsilon^{\star})
s^2_{k+1}(M) \le \max\{ 3 s^2_{k+1}(M),
\frac{3mn\log m}{\ell} \}$. Therefore, 
$$\frac{s_k(U_{1:k}U_{1:k}^\dagger M^{(B)})}{s_{1}((I-U_{1:k}U_{1:k}^\dagger) M^{(B)})} = \omega(1),$$ since
$s^2_k (M) = \omega (\frac{mn\log m}{ \ell })$ and
$\frac{s_k (M)}{s_{k+1} (M)} = \omega(1)$. 

\subsection{Proof of Theorem~\ref{thm:streamingV}}\label{sec:proof-thm5}

We can rewrite
$(A^{(1)})^{\dagger}  W$ as follows:
\begin{eqnarray*}
(A^{(1)})^{\dagger}  W
  &=&\mathbb{E}[(A^{(1)})^{\dagger}]  W
      +((A^{(1)})^{\dagger} - \mathbb{E}[(A^{(1)})^{\dagger}]) 
      W \cr
  &=&V_{1:k}V_{1:k}^\dagger   \mathbb{E}[(A^{(1)})^{\dagger}]  W
      + (I-V_{1:k}V_{1:k}^\dagger)  
      \mathbb{E}[(A^{(1)})^{\dagger}]  W + 
      ((A^{(1)})^{\dagger} - \mathbb{E}[(A^{(1)})^{\dagger}])  W.
\end{eqnarray*}
In the above equation, the columns of 
$(V_{1:k}V_{1:k}^\dagger   \mathbb{E}[(A^{(1)})^{\dagger}]  W)$
have the same space
what we want to recover and the remaining part is noise. Thus, we can easily recover $\hat{V}$ satisfying $\|V_{1:k}^{\dagger}  \hat{V}_{\bot}
\| = o(1)$ when
\begin{equation}\frac{s_k (V_{1:k}V_{1:k}^\dagger  
  \mathbb{E}[(A^{(1)})^{\dagger}]  W)}
{\|(I-V_{1:k}V_{1:k}^\dagger)  
      \mathbb{E}[(A^{(1)})^{\dagger}]  W\|_2+\|((A^{(1)})^{\dagger} - \mathbb{E}[(A^{(1)})^{\dagger}]) 
      W\|_2} = \omega(1).\label{eq:cons}
\end{equation}
Before giving the proof of \eqref{eq:cons} to conclude the proof of Theorem~\ref{thm:streamingV}, we introduce
key lemmas. Lemma~\ref{lem:pf:thm3-1} finds a lower bound for
$s_k (V_{1:k}V_{1:k}^\dagger  
\mathbb{E}[(A^{(1)})^{\dagger}]  W)$
and an upper bound for
$\|(I-V_{1:k}V_{1:k}^\dagger)   \mathbb{E}[(A^{(1)})^{\dagger}] 
W\|_2$
and Lemma~\ref{lem:thm3-2} induces an upper bound for
$\|((A^{(1)})^{\dagger} - \mathbb{E}[(A^{(1)})^{\dagger}])  W\|_2$.
\begin{lemma} When $s_k^2 (M) = \omega(\frac{mn \log m}{\ell})$, $\frac{s_k (M)}{s_{k+1} (M)} = \omega (1)$, and $\frac{\delta^2 \ell s_k^4
  (M)}{ m n^2(k+ \log \ell)} = \omega(1)$, with high probability,
\begin{eqnarray*}
s_k (V_{1:k}V_{1:k}^\dagger  
\mathbb{E}[(A^{(1)})^{\dagger}]  W) & = & \Omega\left( \delta^2 s^2_{k}(M)\sqrt{\frac{\ell}{n}}\right)\quad \mbox{and}\cr
 \|(I-V_{1:k}V_{1:k}^\dagger)   \mathbb{E}[(A^{(1)})^{\dagger}] 
W\|_2& = &o\left( \delta^2 s^2_{k}(M)\sqrt{\frac{\ell}{n}}\right).
\end{eqnarray*}
\label{lem:pf:thm3-1}
\end{lemma}
\noindent{\bf \em Proof:} The proof is given in Section~\ref{sec:proof-lemma-refl}.\QED

\begin{lemma} For given $Q$ and $A^{(B_2)}$, 
$\mathbb{E}[\| ((A^{(1)})^{\dagger} - \mathbb{E}[(A^{(1)})^{\dagger}])  W
\|^2_F ] = O(\delta^2kmn) $.\label{lem:thm3-2}
\end{lemma}
\noindent{\bf \em Proof:} Since every entry of $A^{(1)}$ is randomly sampled with
probability $\delta/4$ and
$W$ and $A^{(1)}$ are independent, for all $1\le i\le n$ and
$1\le j \le k$,
\begin{eqnarray*} 
\mathbb{E}\left[ \big([(A^{(1)} -
\mathbb{E}[A^{(1)}])^\dagger   W]_{ij}\big)^2 \right] &=& \mathbb{E}\left[ \big(\sum_{u=1}^{n}[A^{(1)} -
\mathbb{E}[A^{(1)}]]_{ui}   [W]_{uj}\big)^2 \right]\cr
&\le& \frac{\delta}{4}\|W_j \|^2
= O( \delta^2 m ),
\end{eqnarray*}
where the last equality stems from the trimming process on $A^{(B_2)}$.
Thus, $$\mathbb{E}[\| (A^{(1)} -
\mathbb{E}[A^{(1)}])^\dagger   W \|^2_F ] = O(\delta^2kmn) .$$\QED

\smallskip
\noindent{\bf \em Proof of Theorem~\ref{thm:streamingV}:} When $\frac{\delta^2 \ell s^4_k (M)}{kmn^2} =
\omega(1)$, from Lemma~\ref{lem:pf:thm3-1}, Lemma~\ref{lem:thm3-2}, and the Markov
inequality, $\frac{s_k(V_{1:k}V_{1:k}^\dagger   \hat{V})}{\|(I-V_{1:k}V_{1:k}^\dagger)   \hat{V} \|_2} = \omega(1),$ with high
probability.
Let $\hat{V}= V'\Sigma'(U')^\dagger$ be the singular value
decomposition of $\hat{V}$. Since 
$$\|(I-V_{1:k}V_{1:k}^\dagger)
\hat{V} \|_2 \ge \|(I-V_{1:k}V_{1:k}^\dagger) V' \|_2 s_k(\hat{V}) =
\|(V_{1:k})_\bot^\dagger V' \|_2
s_k(\hat{V})$$
and $s_k(\hat{V}) = \Omega(s_k(V_{1:k}V_{1:k}^\dagger   \hat{V}))$
from the Lidskii ineuality $s_{k+1}(A+B) \ge s_k(A)-S_{k+1}(A)$,
$\frac{s_k(V_{1:k}V_{1:k}^\dagger
  \hat{V})}{\|(I-V_{1:k}V_{1:k}^\dagger)   \hat{V} \|_2} = \omega(1)$
implies $\|(V_{1:k})_\bot^\dagger V' \|_2 = o(1)$. Therefore, with
high probability, 
$$\|V_{1:k}^\dagger (\hat{V})_\bot \|_2 = \sqrt{1- s^2_k(V_{1:k}^\dagger V')}=
\|(V_{1:k})_\bot^\dagger V' \|_2 = o(1) .$$

\subsection{Proof of Lemma~\ref{lem:pf:thm3-1}} \label{sec:proof-lemma-refl}
Since $W = A^{(B_2)}  Q =
\mathbb{E}[A^{(B_2)}]  Q+(A^{(B_2)} -
\mathbb{E}[A^{(B_2)}])  Q$, we find a lower bound for
$s_k (V_{1:k}V_{1:k}^\dagger  
\mathbb{E}[(A^{(1)})^{\dagger}]  W)$
and an upper bound for
$\|(I-V_{1:k}V_{1:k}^\dagger)   \mathbb{E}[(A^{(1)})^{\dagger}] 
W\|_2$ from
\begin{eqnarray}
s_k (V_{1:k}V_{1:k}^\dagger  
\mathbb{E}[(A^{(1)})^{\dagger}]  W) & \ge & s_k (V_{1:k}V_{1:k}^\dagger  
\mathbb{E}[(A^{(1)})^{\dagger}]  \mathbb{E}[A^{(B_2)}]  Q) -\cr
&& \|(\mathbb{E}[A^{(1)}])^\dagger   ((A^{(B_2)} -
\mathbb{E}[A^{(B_2)}])  Q)\|_2 \quad \mbox{and}\cr
 \|(I-V_{1:k}V_{1:k}^\dagger)   \mathbb{E}[(A^{(1)})^{\dagger}] 
W\|_2& \le & \|(I-V_{1:k}V_{1:k}^\dagger)   \mathbb{E}[(A^{(1)})^{\dagger}] 
\mathbb{E}[A^{(B_2)}]  Q\|_2+ \cr
&&\|(\mathbb{E}[A^{(1)}])^\dagger   ((A^{(B_2)} -
\mathbb{E}[A^{(B_2)}])  Q)\|_2.\label{eq:lempf}
\end{eqnarray}

\smallskip
\noindent{\bf \em Key lemmas:} The following lemmas bound each
element of the above inequalities. To show the lemmas, we use Corollary~\ref{cor:1}: $\|(\bar{V}_{1:k}^{1:\ell})^\dagger   Q_\bot\| = o(1)$ with high probability when $\sigma_k^2 (M) = \omega(\frac{mn \log m}{\ell})$, $\frac{s_k (M)}{s_{k+1} (M)} = \omega (1)$, and $\frac{\delta^2 \ell s_k^4
  (M)}{ m n^2\log \ell} = \omega(1)$. 
\begin{lemma}When $s_k^2 (M) = \omega(\frac{mn \log m}{\ell})$, $\frac{s_k (M)}{s_{k+1} (M)} = \omega (1)$, and $\frac{\delta^2 \ell s_k^4
  (M)}{ m n^2\log \ell} = \omega(1)$, with high probability, 
$$ s_k (V_{1:k}V_{1:k}^\dagger  
\mathbb{E}[(A^{(1)})^{\dagger}]  \mathbb{E}[A^{(B_2)}]  Q)
=\Omega\left( \delta^2 s^2_{k}(M)\sqrt{\frac{\ell}{n}}\right). $$\label{lem:9}
\end{lemma}
\noindent{\bf \em Proof:} Since every entry of $A^{(B_2)}$ and $A^{(1)}$ is randomly
sampled with probability $\delta/4$,  we know that $\mathbb{E}[(
A^{(1)})^\dagger] = \frac{\delta}{4} V   \Sigma  
U^{\dagger}$ and $\mathbb{E}[
A^{(B_2)}] = \frac{\delta}{4} U  \Sigma  
(V^{1:\ell})^{\dagger}$. 
Under the conditions of this lemma, from Corollary~\ref{cor:1}
$\|(\bar{V}^{1:\ell})^\dagger   Q_\bot\| = o(1)$ and from Lemma~\ref{lem:cor} $s_k (U_{1:k}^\dagger
M^{(B)}) \ge \sqrt{\frac{\ell}{2n}}s_k(M)$ with high
probability. Let $\bar{R}^{(B)}$ be the $k\times k$ matrix satisfying $V^{1:\ell}_{1:k}=\bar{V}^{1:\ell}\bar{R}^{(B)}$.  Then,
\begin{eqnarray*} 
s_k (  V_{1:k}V_{1:k}^\dagger   \mathbb{E}[(A^{(1)})^\dagger]   (\mathbb{E}[A^{(B_2)}] 
Q) ) &=& \frac{\delta^2}{16} s_k (  V_{1:k}V_{1:k}^\dagger   M^\dagger  M^{(B)} 
Q) )\cr
&=& \frac{\delta^2}{16} s_k (  V_{1:k} \Sigma_{1:k}^{1:k}
    \Sigma_{1:k}^{1:k} (V^{1:\ell})^\dagger Q) )\cr
&\ge& \frac{\delta^2}{16} s_k ( M ) s_k( \Sigma_{1:k}^{1:k}
      (V^{1:\ell})^\dagger Q) )\cr
&=& \frac{\delta^2}{16} s_k ( M ) s_k(
    \Sigma_{1:k}^{1:k}(\bar{R}^{(B)})^\dagger(\bar{V}^{1:\ell})^\dagger
    Q) ) \cr
&\ge&\frac{\delta^2}{16} s_k ( M ) s_k(
    \Sigma_{1:k}^{1:k}(\bar{R}^{(B)})^\dagger) s_k((\bar{V}^{1:\ell})^\dagger
    Q) ) \cr
&=& \Omega\left( \delta^2 s^2_{k}(M)\sqrt{\frac{\ell}{n}}\right), 
\end{eqnarray*}
where the last equality stems from the fact that $s_k(
    \Sigma_{1:k}^{1:k}(\bar{R}^{(B)})^\dagger) = s_k(M^{(B)})$ and $s_k((\bar{V}^{1:\ell})^\dagger
    Q) ) = 1- o(1)$.
 \QED

\begin{lemma}When $s_k^2 (M) = \omega(\frac{mn \log m}{\ell})$, $\frac{s_k (M)}{s_{k+1} (M)} = \omega (1)$, and $\frac{\delta^2 \ell s_k^4
  (M)}{ m n^2\log \ell} = \omega(1)$, with high probability, 
$$ \|(I-V_{1:k}V_{1:k}^\dagger) 
\mathbb{E}[(A^{(1)})^{\dagger}]  \mathbb{E}[A^{(B_2)}]  Q \|_2
=o\left( \delta^2 s^2_{k}(M)\sqrt{\frac{\ell}{n}}\right). $$\label{lem:9-2}
\end{lemma}
\noindent{\bf \em Proof:} Since $\mathbb{E}[(
A^{(1)})^\dagger] = \frac{\delta}{4} V   \Sigma  
U^{\dagger}$ and $\mathbb{E}[
A^{(B_2)}] = \frac{\delta}{4} U  \Sigma  
(V^{1:\ell})^{\dagger}$,
\begin{eqnarray*}
(I-V_{1:k}V_{1:k}^\dagger) 
\mathbb{E}[(A^{(1)})^{\dagger}]  \mathbb{E}[A^{(B_2)}]  Q&=& V_{k+1:n\wedge m}V_{k+1:n\wedge m}^\dagger   \mathbb{E}[(A^{(1)})^\dagger]   \mathbb{E}[A^{(B_2)}]  Q\cr
&=&\frac{\delta^2}{16} V_{k+1:n\wedge m}   \Sigma_{k+1:n\wedge m}^{k+1:n\wedge m} U_{k+1:n\wedge m}^\dagger   M^{(B)}   Q .
\end{eqnarray*}
Under the conditions of this lemman, $s_1 (U_{k+1:n\wedge m}^\dagger
M^{(B)}) = o(\sqrt{\frac{\ell}{n}}\sigma_k(M))$ with high
probability from Lemma~\ref{lem:cor}.
Therefore,
\begin{align*} 
s_1((I-V_{1:k}V_{1:k}^\dagger)  
\mathbb{E}[(A^{(1)})^{\dagger}]  \mathbb{E}[A^{(B_2)}]  Q ) 
=& \frac{\delta^2}{16}s_1 ( V_{k+1:n\wedge m}^\dagger
   \Sigma_{k+1:n\wedge m}^{k+1:n\wedge m}\Sigma_{k+1:n\wedge m}^{k+1:n\wedge m}V_{k+1:n\wedge m}^\dagger Q )
   \cr
\le&\frac{\delta^2}{16}s_{k+1} (M)
     s_1(\Sigma_{k+1:n\wedge m}^{k+1:n\wedge m}V_{k+1:n\wedge m}^\dagger Q ) \cr
\le&\frac{\delta^2}{16}s_{k+1} (M)
     s_1(\Sigma_{k+1:n\wedge m}^{k+1:n\wedge m}V_{k+1:n\wedge m}^\dagger ) \cr
=&o\left( \delta^2 s^2_{k}(M)\sqrt{\frac{\ell}{n}}\right), 
\end{align*}
where the last equality stems from the fact that
$\frac{s_k(M)}{s_{k+1}(M)}=\omega(1)$ and
$s_1(\Sigma_{k+1:n\wedge m}^{k+1:n\wedge m}V_{k+1:n\wedge m}^\dagger ) = s_1 (U_{k+1:n\wedge m}^\dagger
M^{(B)}) = o(s_k(M)\sqrt{\ell/n})$.
 \QED

\begin{lemma} With probability $1-1/\delta$,
$\| (\mathbb{E}[A^{(1)}])^\dagger   ((A^{(B_2)} -
\mathbb{E}[A^{(B_2)}])  Q_{1:k}) \|_2 = O( \sqrt{\delta^2 kmn}).$\label{lem:10}
\end{lemma}
\noindent{\bf \em Proof:} Since entries of $A^{(B_2)}$ are randomly
sampled with probability $\delta/4$ and independent with $Q$, for all $1\le
i\le n$ and $1\le j \le k$,
\begin{align*}
\mathbb{E}\Big[ \big( [(\mathbb{E}[A^{(1)}])^\dagger &  ((A^{(B_2)} -
\mathbb{E}[A^{(B_2)}])  Q)]_{ij} \big)^2 \Big] \cr
=& \mathbb{E}\Big[ \big(\frac{\delta}{4}\sum_{u=1}^{m} \sum_{v=1}^{\ell} [M]_{ui}   [A^{(B_2)} -
\mathbb{E}[A^{(B_2)}]]_{uv}  [Q]_{vj} \big)^2 \Big] \cr
  =&\frac{\delta^2}{16}\sum_{u=1}^{m} \sum_{v=1}^{\ell} 
     [M]_{ui}^2 [Q]_{vj}^2 \mathbb{E}[ ([A^{(B_2)} -
\mathbb{E}[A^{(B_2)}]]_{uv})^2 ] \cr
\le & \frac{\delta^2}{16} \sum_{u=1}^{m} [M]_{ui}^2 \sum_{v=1}^{\ell} 
      [Q]_{vj}^2  \frac{\delta}{4} \le \left(\frac{\delta}{4} \right)^3m .
\end{align*} 
From the above inequality, 
$\mathbb{E}[\| (\mathbb{E}[A^{(1)}])^\dagger   ((A^{(B_2)} -
\mathbb{E}[A^{(B_2)}])  Q) \|^2_F ]= \left( \frac{\delta}{4}
\right)^3 kmn.$ Therefore, by the Markov inequality,
we conclude this proof. \QED

\smallskip
\noindent{\bf \em Proof of Lemma~\ref{lem:pf:thm3-1}:} When $\frac{\delta^2 \ell s_k^4
  (M)}{ m n^2(k+ \log \ell)} = \omega(1)$, $\frac{s_k (V_{1:k}V_{1:k}^\dagger  
\mathbb{E}[(A^{(1)})^{\dagger}]  W)}{\sqrt{\delta^2 kmn}} = \omega(1) $. Inserting Lemma~\ref{lem:9}, Lemma~\ref{lem:9-2}, and
Lemma~\ref{lem:10} into \eqref{eq:lempf}, therefore, we conclude
this proof: 
\begin{eqnarray*}
s_k (V_{1:k}V_{1:k}^\dagger  
\mathbb{E}[(A^{(1)})^{\dagger}]  W) & = & \Omega\left( \delta^2 s_{k}(M^\dagger M)\sqrt{\frac{\ell}{n}}\right)\quad \mbox{and}\cr
 \|(I-V_{1:k}V_{1:k}^\dagger)   \mathbb{E}[(A^{(1)})^{\dagger}] 
W\|_2& = &o\left( \delta^2 s_{k}(M^\dagger M)\sqrt{\frac{\ell}{n}}\right).
\end{eqnarray*}

\subsection{Proof of Theorem~\ref{thm:stmcomp}}

Let $P_{\hat{V}} =
\hat{V}\hat{R}\hat{R}^\dagger \hat{V}^\dagger$ which is an orthogonal
projection matrix onto the linear
span of $\hat{V}$. Then, $\hat{U}\hat{V} =
\frac{4}{\delta} A^{(2)} P_{\hat{V}}$. We can bound $\|
|\hat{U}\hat{V} |^1_0 -M\|_F^2$ using the projection
$P_{\hat{V}}$ as follows:
\begin{eqnarray*}
\| |\frac{4}{\delta} A^{(2)} P_{\hat{V}} |^1_0 -M\|_F^2 
& = & \|
|(M+\frac{4}{\delta}(A^{(2)} - \frac{\delta}{4} M) ) P_{\hat{V}}|^1_0 -M \|_F^2
      \cr
& \le & \|
(M+\frac{4}{\delta}(A^{(2)} - \frac{\delta}{4} M) ) P_{\hat{V}} -M \|_F^2
      \cr
&\le& 2\|MP_{\hat{V}}-M \|_F^2 + 2\|\frac{4}{\delta}(A^{(2)} - \frac{\delta}{4} M)
P_{\hat{V}} \|_F^2 \cr
&\stackrel{(a)}{\le}& 2\|MP_{\hat{V}}-M \|_F^2 + o(mn) \cr
&\le& 4\|U_{1:k}U_{1:k}^\dagger   (MP_{\hat{V}}-M) \|_F^2 + 4\|(I-U_{1:k}U_{1:k}^\dagger)   (MP_{\hat{V}}-M) \|_F^2 + o(mn)\cr
&\stackrel{(b)}{\le}& 4\|U_{1:k}U_{1:k}^\dagger   (MP_{\hat{V}}-M) \|_F^2 + o(mn) \cr 
&=&4\|U_{1:k}\Sigma_{[k]}V_{1:k}^\dagger
                    (P_{\hat{V}}-I) \|_F^2+ o(mn)\cr
&\stackrel{(c)}{=}&o(mn),
\end{eqnarray*}
where $(a)$ stems from Lemma~\ref{lem:x-gam}, $(b)$ uses the fact that
$\|(I-U_{1:k}U_{1:k}^\dagger)M\|_F^2 = o(mn)$, and $(c)$ holds since $\|V^{\dagger}\hat{V}_{\bot} \|_F^2 =
o(1)$ from Theorem~\ref{thm:streamingV}.

\begin{lemma}When $n=\omega(K/\delta)$, with high probability,
$\|\frac{4}{\delta}(A^{(2)} - \frac{\delta}{4} M)
P_{\hat{V}} \|_F^2 =o(mn)$.\label{lem:x-gam}
\end{lemma}
\noindent{\bf \em Proof:} Since entries of $A^{(2)}$ are randomly
sampled with probability $\delta/4$ and independent with $\hat{V}$, for all $1\le
i\le m$ and $1\le j \le n$,
\begin{align*}
\mathbb{E}\Big[\big([(A^{(2)} - \frac{\delta}{4} M)
P_{\hat{V}}]_{ij}\big)^2\Big]
  =&\mathbb{E}\Big[\big(\sum_{v=1}^n[A^{(2)} - \frac{\delta}{4} M]_{iv}
[P_{\hat{V}}]_{vj}\big)^2\Big]
\cr
=&\sum_{v=1}^{n} [P_{\hat{V}}]_{vj}^2 \mathbb{E}[ ([A^{(2)} - \frac{\delta}{4} M]_{iv})^2] \le  \frac{\delta}{4} \sum_{v=1}^{n} 
     [P_{\hat{V}}]_{vj}^2.
\end{align*} 
Since $\sum_{w=1}^{n}\sum_{v=1}^{n}[P_{\hat{V}}]_{vw}^2 = k$, 
from the above inequality, 
$$\mathbb{E}[\|\frac{4}{\delta}(A^{(2)} - \frac{\delta}{4} M)
P_{\hat{V}} \|_F^2] =\left(\frac{4}{\delta}\right)^2  \sum_{i=1}^{m}\sum_{j=1}^n \mathbb{E}\Big[\big([(A^{(2)} - \frac{\delta}{4} M)
P_{\hat{V}}]_{ij}\big)^2\Big] \le \frac{4km}{\delta}.$$
Therefore, by the Markov inequality,
we conclude this proof. \QED

\end{document}